\documentstyle[12pt]{article}
\input amssym.def
\input amssym
\title{{\bf Rank one lattice type vertex operator algebras and their
automorphism groups}
\footnotetext{The first author is supported by NSF grant DMS-9700923 and a
research grant from the Committee on Research, UC Santa Cruz.}\footnotetext{
The second author is supported by NSF grant DMS-9623038 and the
University of Michigan Faculty Recognition 
Grant (1993--96).}}
\author{Chongying Dong and  Robert L.~Griess Jr.}
\date{}

\makeatletter

\begin{document}

\newtheorem{thm}{Theorem}[section]
\newtheorem{prop}[thm]{Proposition}
\newtheorem{lem}[thm]{Lemma}
\newtheorem{rem}[thm]{Remark}
\newtheorem{cor}[thm]{Corollary}
\newtheorem{conj}[thm]{Conjecture}
\newtheorem{de}[thm]{Definition}
\newtheorem{notation}[thm]{Notation}
\pagestyle{plain}
\maketitle

\def\vv{ {\hbox {\bf 1} }}  
\def \Z{\Bbb Z}
\def \F{\Bbb F}
\def \C{\Bbb C}
\def \R{\Bbb R}
\def \Q{\Bbb Q}
\def \N{\Bbb N}
\def \D{{\cal D}}
\def \wt{{\rm wt}}
\def \tr{{\rm tr}}
\def \sp{{\rm span}}
\def \Res{{\rm Res}}
\def \Res{{\rm QRes}}
\def \End{{\rm End}}
\def \E{{\rm End}}
\def \Ind {{\rm Ind}}
\def \Irr {{\rm Irr}}
\def \Aut{{\rm Aut}}
\def \Hom{{\rm Hom}}
\def \mod{{\rm mod}}
\def \ann{{\rm Ann}}
\def \<{\langle} 
\def \>{\rangle} 
\def \t{\tau }
\def \a{\alpha }
\def \e{\epsilon }
\def \l{\lambda }
\def \L{\Lambda }
\def \g{\frak g}
\def \b{\beta }
\def \om{\omega }
\def \o{\omega }
\def \k{\kappa}
\def \c{\chi}
\def \ch{\chi}
\def \cg{\chi_g}
\def \ag{\alpha_g}
\def \ah{\alpha_h}
\def \ph{\psi_h}
\def \pf{{\bf Proof. }}
\def \voa{{vertex operator algebra\ }}
\def \svoa{{super vertex operator algebra\ }}
\def \qed{\mbox{$\square$} \vskip .6cm }
\def \lc{L_C}
\def \tlc{\widetilde{L}_C}
\def \tv{\widetilde{V}_L}
\def \vlc{V_{L_C}}
\def\tvlc{\widetilde{V}_{L_C}}
\def\vtlc{V_{\widetilde{L}_C}}
\def\tvtlc{\widetilde{V}_{\widetilde{L}_C}}
\def\ha{\frac{1}{2}}
\def\se{\frac{1}{16}}

\def\Ve{V^{0}}
\def\M{\rm \ {\sl  M}\llap{{\sl I\kern.80em}}\ }
\def\xx{\em}
\def\s{\sigma}
\def\a{\alpha} 
\def\t{\tau}
\def\b{\beta} 

\def\xap{x_\a^+}
\def\xam{x_\a^-}

\def\la{\langle}
\def\ra{\rangle}
\def\rtar{\rightarrow}
\def\mt{\mapsto}


{\bf Abstract.} Let $L$ be a positive definite even lattice of rank one
and $V_L^+$ be the fixed points of the lattice VOA $V_L$ associated to 
$L$ under an  automorphism of $V_L$ lifting the $-1$ isometry of
$L.$ A set of 
generators and the full automorphism group of $V_L^+$ are determined.  

\section{Introduction}

The classification of rational vertex operator algebras is definitely
one of the most important steps toward classification of
rational conformal field theory. If the central charge $c<1$ there is a
complete classification of the unitary representations of the 
Virasoro algebra [GKO], [FQS]. It turns out 
the associated Virasoro vertex operator algebras
are the only rational VOAs  such that the action of the Virasoro algebra
is unitary [DMZ].  
Moreover the vertex operator algebra associated 
the irreducible highest weight representations of the Virasoro
algebra with central charge $c$ is rational if and only if 
$c=1-6\frac{(p-q)^2}{pq}$ for two coprime positive integers $p,q>1$
[W].   

It is natural to consider classifying the rational vertex operator algebras
with $c=1$. A lot of analysis  has been done in the physics literature
in this direction [G], [H], [DVVV], [K]. In particular, the partition
functions of rational conformal field theories  with $c=1$ have been
classified in [K] as proposed in [G].  

If we hope to classify rational VOAs with $c=1$, it is necessary 
to have a better understanding of the known rank 1 VOAs.  
This is exactly our goal in this paper. We study the lattice type 
vertex operator algebras with $c=1$ and determine their automorphism
groups.  

The basic definitions for lattice VOAs are given in [B] and [FLM].  This VOA is
denoted
$V_L$, where $L$ is an even integral lattice of finite rank.  There is a
family
of  automorphisms $\theta$ of $V_L$, unique up to conjugacy, which
interchanges
the 1-spaces spanned by all $e^\lambda$ and $e^{-\lambda}$ [DGH], Appendix
D;
these are the ``lifts of minus the identity".  We may arrange for such a
lift
to satisfy  $e^\lambda \mapsto e^{-\lambda}$, for all $\lambda$.   Call
this
automorphism the ``standard lift of minus the identity", and denote its
fixed
points on $V_L$ by $V_L^+$. We call this fixed point VOA a  {\sl
symmetrized
lattice VOA}. 
For all rank one lattices, 
we determine the automorphism group and give sets of generators for the
symmetrized lattice VOA.  The analogous problem for lattice VOAs has been
solved in [DN]. We use the term {\sl lattice type VOA} for a lattice
VOA or symmetrized lattice VOA.

\bigskip

We may  take any subgroup $F$ of $Aut(V)$ and consider the subVOA of fixed
points, $V^F$ (this VOA also has rank 1).  
The group $N_G(F)$ acts as automorphisms of it and the
kernel of
this action contains $F$.   In some cases, $V^F$ turns out to be a lattice
VOA or a symmetrized lattice VOA.   
This is the case where $V$ is a
lattice VOA
and  $F$ is an extension $A:B$, where $A$ is a subgroup of 
the natural maximal  torus associated to the lattice 
and $B$ is a group of order 2 generated 
by an element
of the
torus normalizer acting as $-1$ on it.   The rank 1 lattices of Section 2
occur
this way within $V_{A_1}$,    the lattice VOA for the $A_1$-lattice.  The
advantage of this viewpoint is that we can see automorphisms of the subVOA
which
are not obvious from the definition.  For example, the exceptional case
$n=4$ of
Theorem 3.1 can be anticipated as follows.  If we take $F$ to be a four
group
in  $G:=Aut(V_{A_1}) \cong PSL(2,C)$, then $F=C_G(F)$ and $N_G(F) \cong
Sym_4$
and so we get an action of $N_G(F)/F  \cong Sym_3$ on $V^F \cong V_L$,
where $L
\cong 2L_{A_1}$. This action in fact is faithful; see
details in Section 4. 
The automorphism results of this paper may be combined with
[DN] to say that  every finite subgroup $F$ of $Aut(V_{A_1})$
which is of ``AD" type (cyclic or dihedral), $N_G(F)/F$ induces $Aut(V^F)$
on
$V^F$, a pleasant situation.   We are currently 
studying the ``E type" examples,
i.e., $F
\cong Alt_4, Sym_4, Alt_5$.  

\bigskip

The paper is organized as follows: In Section 2 we review the construction
of lattice vertex operator algebras $V_L$ associated to rank one lattice 
from [B] and [FLM]. We decompose $V_L$ into a direct sum of irreducible
modules for the Virasoro algebra with central charge 1. Moreover,
in the case that $L$ is the root lattice of type $A_1$ we prove that
the $\frak g$-invariants of $V_L^{\frak g}$ is isomorphic to the
highest weight irreducible module $L(1,0)$ for the Virasoro algebra
where $\frak g \cong sl(2,\C)$ is the Lie algebra 
$(V_1,\ 1^{st} \hbox{ binary composition})$.  We also apply 
a result in [DLM] to decompose $V_L$ into irreducible modules
for the pair $(\hbox{Virasoro algebra}, \frak g).$ The main result in
this section is the determination of a set of generators for
$V_L^+$ for an arbitrary even positive definite lattice of rank 1.
Several results in Section 2 have already appeared in the physics paper
[DVVV].  

Sections 3 and 4 are devoted to the automorphism groups of 
$V_L^+.$ The most difficult case is when $L$ is generated by
an element of square length 8. As discussed before, we identify
$V_L^+$ in this case with $V_{A_1}^{F}$ whose automorphism
group is relatively easy to determine by using
 knowledge of Lie algebra $sl(2,\C).$ Section 5 is an appendix 
on commutative nonassociative 
algebras of dimension $n-1$ with $Sym_n$ as automorphism
group. This result is used in Section 4 to determine 
$Aut (V_{L_8}^+).$ 

\section{Vertex operator algebras $V_L^+$ }\label{one}
\setcounter{equation}{0}

In this section we first 
briefly review the construction of rank one lattice VOA
from [FLM]. Let $L$ be a rank one even positive definite lattice,
${\frak h}=L\otimes_{\Bbb Z}{\Bbb C},$ and 
$\hat{\frak h}_{\Bbb Z}$  the corresponding Heisenberg algebra;
the bilinear form on $L$ or $\frak h$ is denoted $\< \cdot , \cdot \>$.
Denote by  $M(1)$  the associated irreducible induced module for 
$\hat{\frak h}_{\Bbb Z}$ 
such that the canonical central element of $\hat{\frak h}_{\Bbb Z}$ acts as 1;
it is identified with the ``polynomial factor" in $V_L$. 
Let $\C[L]$ be the group algebra of $L$ with a basis 
$e^{\a}$ for $\alpha\in L.$ Let $\beta\in \frak h$ such that 
$\<\beta,\beta\>=1.$ It was proved in~[B] and [FLM] that there is a linear map 
$$\begin{array}{lcr}
V_{L}&\to& (\mbox{End}\,V_{L})[[z,z^{-1}]],\hspace*{3.6 cm} \\
v&\mapsto& Y(v,z)=\displaystyle{\sum_{n\in\Z}v_nz^{-n-1}\ \ \ (v_n\in
\mbox{End}\,V_{L})}
\end{array}$$ such that
$V_{L}=(V_{L},Y,{\bf 1},\omega)$ is a simple vertex operator algebra where
${\bf 1}=1\otimes e^0$ and $\omega=\frac{1}{2}\beta(-1)^2.$  

Let $\theta$ be an  automorphism of $L$ such that 
$\theta(\a)=-\a$; see Appendix D of [DGH]. We define an automorphism of $V_L$,
denoted again by $\theta$, such that $\theta (u\otimes e^{\alpha})
=\theta(u)\otimes e^{-\alpha}$ for $u\in M(1)$ and $\a\in L$.
Here the action of $\theta$ on $M(1)$ is given by 
$\theta(\beta(n_1)\cdots
\b(n_k))=(-1)^k\b(n_1)\cdots\b(n_k)$. The 
$\theta$-invariants $V_L^+$ of $V_L$ form a simple vertex operator 
subalgebra and the $(-1)$-eigenspace $V_L^-$ is an irreducible $V_L^+$-module
(see Theorem~2 of [DM2]). Clearly $V_L=V_L^+ \oplus V_L^-$.

\begin{rem}{\rm Clearly,  there exists an even positive integer $n$ such that 
$L \cong L_n :=  \Z\sqrt{n}\beta.$ If $n=2$, 
$V_L$ affords  the fundamental representation 
for the affine algebra $A_1^{(1)}.$}
\end{rem}

It is well-known that both $M(1)$ and $V_L$ are unitary representations for 
the Heisenberg algebra and for the Virasoro algebra (see e.g. [KR]).
We next
discuss the decomposition of $M(1)$ and $V_L$ into 
irreducible modules for the Virasoro algebra. Let $W$ be a module for the Virasoro algebra
with central charge $c$ such that $W=\oplus_{n\in\C}W_n$ where
$W_n$ is the eigenspace for $L(0)$ with eigenvalue $n$ and 
is finite-dimensional. We define the $q$-graded dimension of $W$
as $$\dim_qW=q^{-c/24}\sum_{n\in\C}(\dim W_n)q^n.$$

Denote by $L(c,h)$ the unique irreducible highest weight module for
the Virasoro algebra with central charge $c\in\C$ and highest weight
$h\in\C.$ Then 
$$\dim_qL(1,h)=\left\{\begin{array}{ll}
{1\over \eta(q)}(q^{n^2/4}-q^{(n+2)^2/4}), & {\rm if} \ 
h=\frac{1}{4}n^2,n\in\Z\\
\frac{1}{\eta(q)}q^h, & {\rm otherwise.}
\end{array}\right.$$
(cf. [KR]) where
$$\eta(q)=q^{1/24}\prod_{n=1}^{\infty}(1-q^n).$$ 

For a fixed even positive integer $n,$
$$V_{L_n}=\bigoplus_{m\in \Z}M(1)\otimes e^{m\sqrt{n}\beta}$$
is a decomposition of $V_{L_n}$ into submodules for the Virasoro
algebra. It is clear from the construction that 
$\dim_qM(1)\otimes e^{m\sqrt{n}\beta}={q^{m^2n/2}\over \eta(q)}$ (cf.
[FLM]).  Irreducibility of the module $M(1)\otimes e^{m\sqrt{n}\beta}$
for the Virasoro algebra 
depends on whether 
$2n$ is a perfect square (see e.g. [KR]).
If $2n$ is not a perfect square and $m\ne 0$ then
 $M(1)\otimes e^{m\sqrt{n}\beta}$
is irreducible and isomorphic to $L(1,\frac{m^2n}{2})$ as
they have the same graded dimension. 
Now we assume that $2n=4k^2$ for some nonnegative integer $k.$
Then $\frac{m^2n}{2}={(mk)^2}.$ We have
$${q^{m^2n/2}\over \eta(q)}=\sum_{p\geq 0} 
{1\over \eta(q)}(q^{(mk+p)^2}-q^{(mk+p+1)^2})$$
and accordingly,
$$M(1)\otimes e^{m\sqrt{n}\beta}=\bigoplus_{p\geq 0} L(1,(mk+p)^2).$$

To summarize, we have
\begin{prop}\label{p1} As modules for the Virasoro algebra
$$M(1)=\bigoplus_{p\geq 0}L(1,p^2)$$
$$V_{L_n}=\left(\bigoplus_{p\geq 0}L(1,p^2)\right)\bigoplus \left(\bigoplus_{m>0}2L(1,\frac{m^2n}{2})\right)$$
if $2n$ is not a perfect square, and
$$V_{L_n}=\bigoplus_{m\geq 0}\bigoplus_{p=0}^{k-1}(2m+1)L(1,(mk+p)^2)$$
if $2n=4k^2$ for some nonnegative integer $k.$
\end{prop}

Note that the lattice $L_n$ with $n=2k^2$ is a sublattice of $L_2.$
We now restrict ourselves to the lattice $L_2.$ From Proposition
\ref{p1},  we have the decomposition
\begin{equation}\label{2.1}
V_{L_2}=\bigoplus_{m\geq 0}(2m+1)L(1,m^2).
\end{equation}
Note that the weight one subspace $(V_{L_2})_1$ of $V_{L_2}$ forms
a Lie algebra $\frak g$  isomorphic to $sl(2,\C).$
Thus, $\frak g$ acts on $V_{L_2}$ via $v_0$ for $v\in (V_{L_2})_1$
where $v_0 \in End(V)$ is defined by $Y(v,z)=\sum_{n\in\Z}v_nz^{-n-1}.$ 
It is clear that $\g$ acts  on $V_{L_2}$ as 
derivations in the sense that for $d\in \frak g$ and $v\in V_{L_2}$ 
$$[d,Y(v,z)]=Y(dv,z).$$ 
The $\frak g$-invariants $V_{L_2}^{\g}=\{v\in V_{L_2}|\g\cdot v=0\}$
form a simple vertex operator algebra (see [DLM]). 

Let $W_m$ be the unique highest weight module for $\g$ with
highest weight $m$; it has dimension $m+1$.  Let 
$V_{L_2}^{W_{m}}$ be 
the sum of irreducible $\g$-submodule
of $V_{L_2}$ isomorphic to $W_m,$ and $(V_{L_2})_{W_{m}}$ 
the space of highest weight vectors in $V_{L_2}^{W_{m}}.$
Then by Theorem
2.8 and Corollary 2.9 of [DLM] or a similar result in [BT], 
as ($V_{L_2}^{\g}$, $\frak g$)-modules $V_{L_2}$ has decomposition
\begin{equation}\label{2.2}
V_{L_2}=\bigoplus_{m\geq 0}V_{L_2}^{W_{2m}}=\bigoplus_{m\geq 0}(V_{L_2})_{W_{2m}}\otimes W_{2m}
\end{equation}
and $(V_{L_2})_{W_{2m}}$ is an irreducible module for $V_{L_2}^{\g}.$ 
Moreover, $(V_{L_2})_{W_{2k}}$ and $(V_{L_2})_{W_{2m}}$ are isomorphic if and
only if $k=m.$  

\begin{rem}{\rm  There was a mistake in Theorem 2.8 of [DLM]. 
The correct statement should be the following
with the same notation as in [DLM]: There exists a set $Q$ of weights
containing the dominant weights in the root lattice of ${\frak g}$ such that
$V^{\l}$ is non-zero if and only if  $\l\in Q,$ and 
$V_{\l}$ is an irreducible $V^{\frak g}$-module. Moreover, $V_{\l}$ 
and $V_{\mu}$ are isomorphic $V^{\frak g }$-modules if, and only, if
$\l=\mu.$ In particular,
$V^{\frak g}$ is a simple vertex operator algebra. In this theorem,
$\g$ is a simple Lie algebra.}
\end{rem}

We immediately have the
following corollary of (\ref{2.1})
\begin{cor}\label{c2.1} $V_{L_2}^{\g}$ and $L(1,0)$  are isomorphic VOAs
and $(V_{L_2})_{W_{2m}}$ is isomorphic to $L(1,m^2)$ as modules.
\end{cor}

Now we are in a position to determine a set of generators for 
vertex operator algebra $M(1)^+=M(1)\cap V_{L_2}^+$ by using 
(\ref{2.1}) and Corollary \ref{c2.1}. For this purpose
we need a detailed description of the module $W_m.$ 
Let $\{x_{\a},\a,x_{-\a}\}$ be the Chevalley basis of $\frak g$ with bracket
$[\a,x_{\pm\a}]=\pm2x_{\pm\a}$  and $[x_{\a},x_{-\a}]=\a.$  
For any weight module $W$ for $\frak g$,  
set $W^i:=\{w\in W|\alpha w=iw\},$ the weight space
of $W$ with weight $i$ for $i\in \C.$ Then 
$W_m=W_{m}^{-m}\oplus W_{m}^{-m+2}\oplus\cdots W_{m}^{m-2}\oplus W_{m}^{m}.$
Since $W_m$ is a unitary representation of $SU(2)$ there is a hermitian
form $(\cdot,\cdot)$ on $W_m$ such that 
$W_m$ has an orthonormal basis 
$\{w_{m,i}|i=-m,-m+2,...,m-2,m\}$ and 
$w_{m,i}$ is a weight vector for $\frak g$ with weight $i.$

Let $m\geq n.$ The following tensor product decomposition
is well-known: 
\begin{equation}\label{2.3}
W_m\otimes W_n=W_{m-n}\oplus W_{m-n+2}\oplus\cdots
\oplus W_{m+n-2}\oplus W_{m+n}.
\end{equation}
It is clear that there is a hermitian form $(\cdot,\cdot)$
on $W_m\otimes W_n$ such that  $\{w_{m,s}\otimes w_{n,t}\}$ forms an orthogonal
basis and the union of the standard bases of $W_i$ also forms an
orthogonal basis. 

\begin{lem}\label{l2.1} If both $m$ and $n$ are even, $m \ge n$,  and 
$$w_{m,0}\otimes w_{n,0}=\sum_{m-n\leq i\leq m+n, i\in 2\Z}c^{mni}_{000}
w_{i,0}$$
then $c^{mni}_{000}=(w_{m,0}\otimes w_{n,0},
w_{i,0})$ is nonzero if and only if ${i\over 2}+{m\over 2}+{n\over 2}$
is even. 
\end{lem}

\pf The constants $c^{mni}_{000}$ are essentially the $SU(2)$ 
Clebsh-Gordan coefficients. For example, one can find this result in 
Section 3.7 of [E].
\qed

Now consider the $\frak g$-submodule $U_{2m}$ inside $V_{L_2}$ generated by
$e^{m\sqrt{2}\beta}$ which is a highest weight vector for
$\g$ with highest weight $2m.$ Then $U_{2m}$ is isomorphic to
$W_{2m}$ as $\frak g$-modules. 
Denote the corresponding basis by $\{u_{2m,j}|j=-2m,-2m+2,...,2m-2,2m\}.$
For convenience set $u^{m}=u_{2m,0}.$ Then $u^m$ necessarily lies in 
$M(1) \cap (V_{L_2})_{m^2}.$ 
Since $e^{m\sqrt{2}\beta}$ is a highest weight vector for the
Virasoro algebra, and the actions of Virasoro algebra and 
of $\frak g$ commute on $V_{L_2}$ we see that $U_{2m}$ is a space
of highest weight vectors for the Virasoro algebra with highest weight
$m^2.$ Thus $u^m$ generates the irreducible highest weight module
isomorphic to $L(1,m^2)$ for the Virasoro algebra in $M(1).$ We will
also use the notation $L(1,m^2)$ for the submodule generated by $u^m$.  

\begin{lem}\label{l2.2} For $m \ge n\geq 0$,  the
submodule of $M(1)$ generated by the coefficients $Y(u^m,z)u^n$
for the Virasoro algebra is exactly 
$L(1,(m-n)^2)\oplus L(1,(m-n+2)^2)\oplus\cdots\oplus L(1,(m+n-2)^2)\oplus
L(1,(m+n)^2).$
\end{lem}

\pf Using (\ref{2.2})-(\ref{2.3}) and Lemma \ref{l2.1} shows that the
submodule of $M(1)$ generated by the coefficients $Y(u^m,z)u^n$
for the Virasoro algebra is contained in  
$L(1,(m-n)^2)\oplus L(1,(m-n+2)^2)\oplus\cdots\oplus L(1,(m+n-2)^2)\oplus
L(1,(m+n)^2).$
Since $L(1,p^2)$ and $L(1,q^2)$ are inequivelent irreducible modules
for the Virasoro algebra if $0\leq p<q$ it is enough to show that 
for any nonnegative integer $i$ congruent to  $m+n$ modulo 2
with $m-n\leq i\leq m+n$ there exists some $k\in\Z$ such that 
the projection of $(u^m)_k(u^n)$ into $L(1,i^2)$ is nonzero. 

This was essentially proved in Lemma 3.1 and Theorem 2 of [DM2].
Although the arguments in [DM2] were for finite groups,  they also 
work for Lie algebras. Thus, for any $i$ with $m-n\leq i\leq m+n$ and
$m+n+i$ even, there exists $p\in\Z$ such that the span $U_{m,n}^p$ 
of $(u_{2m,s})_pu_{2n,t}$ 
for $s=-2m,-2m+2,...,2m-2,2m$ and $t=-2n,-2n+2,...,2n-2,2n$ 
contains a copy $W$ of $W_{2i}.$ By Lemma \ref{l2.1}, the
projection of $(u^{m})_p(u^{n})$ into $W$ is nonzero.
Clearly the image of the projection is in $L(1,i^2).$ 
\qed 

Now we define a linear isomorphism $\theta'$ on $M(1)$ such that
$\theta'$ acts on $L(1,m^2)$ as $(-1)^m.$ Then Lemma \ref{l2.2}
implies that $\theta'$ in fact is a  VOA automorphism of $M(1)$ 
of order 2. 
Since $\theta'$ acts on weight one space as $-1$
and $M(1)$ is generated by $M(1)_1$ we see immediately that
$\theta'$ is exactly the the restriction of
the automorphism $\theta$ of $V_L$.

\begin{thm}\label{t2.1} (1) We have decompositions
$$M(1)^+=\bigoplus_{m\geq 0}L(1,4m^2),\ \ M(1)^-=\bigoplus_{m\geq 0}L(1,(2m+1)^2).$$

(2) $M(1)^+$ is generated by $u^{n^2}$ and $\omega$ for any even positive
integer $n.$ 
\end{thm}

The theorem is an immediate consequence of Lemma \ref{l2.2} and Proposition
\ref{p1}.

\begin{cor} $L(1,0)$ is the only proper vertex operator subalgebra
of $M(1)^+.$
\end{cor}

\pf Let $W$ be a subalgebra of $M(1)^+$, which necessarily contain
$L(1,0)$ as $\omega$ is an invariant for any automorphism and
$L(1,0)$ is generated by $\omega.$ 
If $W$ is strictly greater than $L(1,0)$ then
$W$ contains at least one $L(1,n^2)$ for some $n>0.$ By Theorem 
\ref{t2.1} $L(1,n^2)$ generates $M(1)^+.$ This shows $W$ is
equal to $M(1)^+.$ \qed

Now we turn our attention to $V_{L_n}^+.$ Consider $V_{\Z}:=V_{\Z\beta}$
which is a completely reducible module for $M(1):$
$$V_{\Z}=\bigoplus_{m\in \Z}^{\infty}M(1)\otimes e^{m\beta}$$
where each $M(1)\otimes e^{m\beta}$ is an irreducible $M(1)$-module.
We can define $\theta$ on $V_{\Z}$ as before so that
$\theta e^{m\beta}=e^{-m\beta}$ and the action on $M(1)$ is the same. 
For each $\theta$-stable subspace $W$ of $V_{\Z}$ we define
$W^{\pm}$ to be the eigenspaces with eigenvalues $\pm1.$ Then
$$(M(1)\otimes e^{\beta m}+M(1)\otimes e^{-m\beta})^+
= M(1)^+\otimes (e^{m\beta}+e^{-m\beta})+M(1)^-
\otimes (e^{m\beta}-e^{-m\beta}).$$

Following [DM1] we define a new $M(1)$-module 
$\theta \circ M(1)\otimes e^{m\beta}$ which has the same underlying 
space $M(1)\otimes e^{m\beta}$ with vertex operators
$Y(\theta u,z)$ for $u\in M(1).$ It is easy to see that
$\theta\circ M(1)\otimes e^{m\beta}$ and $M(1)\otimes e^{-m\beta}$ are nonisomorphic
$M(1)$-modules. It follows from Theorem 6.1 of [DM2] that
 $M(1)\otimes e^{m\beta},$ 
$M(1)\otimes e^{-m\beta}$ are isomorphic and irreducible $M(1)^+$-modules. 
In particular, they are isomorphic $L(1,0)$-modules.

{}From the discussion above we see that
\begin{equation}\label{add}
V_{L_{2n}}^+=M(1)^+\bigoplus \left(\bigoplus_{m>1}(V_{L_{2n}}^{+})^m\right)
\end{equation}
where 
\begin{eqnarray*}
& &(V_{L_{2n}}^+)^m=(M(1)\otimes e^{m\sqrt{2n}\beta}+M(1)\otimes e^{-m\sqrt{2n}\beta})^+\\
& &\ \ \ \ \ =M(1)^+\otimes (e^{m\sqrt{2n}\beta}+e^{-m\sqrt{2n}\beta})+M(1)^-\otimes (e^{m\sqrt{2n}\beta}-e^{-m\sqrt{2n}\beta}).
\end{eqnarray*}

Set $e^n:=e^{\sqrt{2n}\beta}+e^{-\sqrt{2n}\beta}.$ 
\begin{thm}\label{t2.2} $V_{L_{2n}}^+$ is generated by $u^4,$ $e^n$ and $\o.$
\end{thm}

\pf  We prove by induction on $m$ that $(V_{L_{2n}}^+)^m$ can be generated 
from  $u^4,$ $e^n$ and $\o.$ The case $m=0$ was proved in Theorem 
\ref{t2.1}. Now we assume $m>0.$  By induction
assumption, $e^{(m-1)\sqrt{2n}\beta}+e^{-(m-1)\sqrt{2n}\beta}$
can be generated.  Then the coefficient 
of $z^{2n(m-1)}$ in 
$$Y(e^n,z)(e^{(m-1)\sqrt{2n}\beta}+e^{-(m-1)\sqrt{2n}\beta})$$
is $e^{m\sqrt{2n}\beta}+e^{-m\sqrt{2n}\beta}+v$
where $v$ lies in $\sum_{k<m}(V_{L_{2n}}^+)^k.$ Thus $e^{m\sqrt{2n}\beta}+e^{-m\sqrt{2n}\beta}$ can be generated. Since $(V_{L_{2n}}^+)^m$ is an 
irreducible
$M(1)^+$-module, the result follows.
\qed

The following Lemma will be needed in the proof of Theorem \ref{t3.1}.
\begin{lem}\label{l2.10} Assume that $4n$ is not a perfect square. 
Then $(V_{L_{2n}}^+)^m$ is an irreducible $L(1,0)$-module isomorphic to 
$L(1,nm^2).$ Moreover, we have the following decomposition
$$(V_{L_{2n}})^+ \cong 
\left(\bigoplus_{p\geq 0}L(1,4p^2)\right)\bigoplus \left(\bigoplus_{m>0}L(1,m^2n)\right)$$
\end{lem}

\pf By the discussion before Proposition \ref{p1} we see that 
both $M(1)\otimes e^{m\sqrt{2n}\beta}$ and $M(1)\otimes e^{-m\sqrt{2n}\beta}$
are irreducible $L(1,0)$-module isomorphic to $L(1,m^2n).$ Note that
$\theta: M(1)\otimes e^{m\sqrt{2n}\beta}\to M(1)\otimes e^{-m\sqrt{2n}\beta}$
is an $L(1,0)$-module isomorphism. Thus 
$$(V_{L_{2n}}^+)^m=\{v+\theta(v)|v\in M(1)\otimes e^{m\sqrt{2n}\beta}\}$$
is also isomorphic to  $L(1,m^2n).$ The decomposition of  $(V_{L_{2n}})^+$
into irreducible $L(1,0)$-modules now follows from Theorem \ref{t2.1} and 
equation (\ref{add}) immediately. \qed

\section{Automorphism groups of $V_{L_{2n}}^+$}
\setcounter{equation}{0}

The main result of this paper is the following.
\begin{thm}\label{t3.1} The full automorphism group of $V_{L_{2n}}^+$
is isomorphic to $\Z_2$ if $n>4$ or $n=3, n=2;$ for $n=1$, it is a semidirect product 
 $T \<\theta\>$, where 
$T$ is a rank 1 torus and  
$\theta$ acts on $T$ as $-1$; 
for $n=4$, it is isomorphic to  $Sym_3$.   
\end{thm}

We now discuss the proof. All cases  except $n=4$ 
will be treated here, and the case $n=4$ will be established in Section 4.  

By Theorem \ref{t2.2} $V_{L_{2n}}^+$ is generated by $u^4$ and $e^n$ as
modules for the Virasoro algebra. Since any automorphism $\sigma$ acts on 
$L(1,0)$ trivially,  we need to determine only the images 
$\sigma(u^4)$ and $\sigma(e^n).$ 

First we discuss the cases $n>4, n=2$ or $n=3.$ 
In the cases $n>4,$ the space of highest weight vectors for the Virasoro 
algebra
in $(V_{L_{2n}}^+)_4$ is one-dimensional and spanned by $u^4.$ In the cases
$n=2$ or $3,$ since $4n$ is not a perfect square we see from Lemma \ref{l2.10}
that $u^4$ is the unique  highest weight vectors for the Virasoro 
algebra
in $(V_{L_{2n}}^+)_4$ up to a scalar multiple. 
Thus
$\sigma(u^4)=cu^4$ for some constant $c.$  By Lemma \ref{l2.2}
we see that $c^2=c$, whence $c=1.$ This shows that the restriction of $\sigma$
to $M(1)^+$ is the identity. 

In order to determine $\sigma (e^n)$ we recall from
Section 2  that $e^n\in (V_{L_{2n}}^+)^1$ and $(V_{L_{2n}}^+)^1$
is an irreducible $M(1)^+$-module. Since $(V_{L_{2n}}^+)^p$
and $(V_{L_{2n}}^+)^q$ are inequivelent irreducible $M(1)^+$-modules
(the lowest weight vectors in these two spaces 
have different weights) if $p<q$ and the space of the lowest weight 
vectors in  $(V_{L_{2n}}^+)^1$ is one dimensional and spanned by
$e^n$ we see immediately that $\sigma(e^n)=de^n$  for a scalar $d.$ 
Note that the coefficient of $z^{-2n}$ in $Y(e^n,z)e^n$ is $2{\bf 1}.$
Thus $d=\pm 1.$ If $d=1$, $\sigma$ is the identity. 
If $d=-1$ the $\sigma$ is the the restriction of automorphism 
$(-1)^{\beta\over\sqrt{2n}}$ of
$V_{L_{2n}}$ to $V_{L_{2n}}^+.$ Here for any $x\in \C\beta$ we define
the operator $e^{\pi ix}$ on $V_{L_{2n}}$ by
 $$ e^{\pi ix}(u\otimes e^{\beta})=e^{\pi i\<x,\beta\>}(u\otimes e^{\beta}).$$
Then  $e^{\pi ix}$ is an automorphism of $V_{L_{2n}}$ (see [DM1]). 
This finishes the proof of Theorem in the case that $n>2.$

Now we deal with the case $n=1.$ In this case $V_{L_2}$ is the
fundamental module for the affine Lie algebra $A_1^{(1)}.$ 
We define 
$\xap := x_\a + x_{- \a}$ and $\xam := x_a - x_{- \a}$.  Then we
have 
$[\xap, \xam]=-2\a$ and for any scalar $c$, $[\a, c \xap]=2c\xam$ and 
$[\a, c \xam]=2c\xap$. Consider $\xap$ as a generator for 
a Cartan subalgebra. Then $\xam$ and $\a$
span the sum of its root spaces. The automorphism can 
be realized as $e^{i\pi (\xap)_0/2}.$ Thus $V_{L_2}^+$ is
isomorphic to the lattice VOA $V_{L_8}.$ So our problem
reduces to determining the automorphism group of $V_{L_8}.$ 

The automorphism group for arbitrary lattice VOA has been determined
in [DN].    But in the
rank one 
case, it is very easy to write down the automorphism group 
directly, and we do so as follows.

For $c\in \C^{\times}$ we define a linear map $\sigma_c$ on 
$V_{L_8}$ by $\sigma_c(u\otimes e^{n\a})=c^n(u\otimes e^{n\a})$ for
$u\in M(1)$ and $n\in\Z.$ It is a straightforward verification 
to show that $\sigma_c$ is indeed an automorphism of $V_{L_8}.$
Since $\theta$ inverts each $\sigma_c$ under conjugation, 
we get 
a subgroup of $\Aut (V_{L_8})$ isomorphic to the semidirect
product  
$\C^\times \colon 2$.  

Let $\sigma$ be an automorphism of $V_{L_8}.$ Then $\sigma\beta(-1)$
=$c\beta(-1)$ as $(V_{L_8})_1$ is one-dimensional and spanned by
$\beta(-1).$ But $\o=\frac{1}{2}\beta(-1)^2$ is an invariant we
see that $c=\pm 1.$ If $c=-1$,  we can multiply $\sigma$ by $\theta.$
So we can assume that $c=1.$ As a result we
have $\sigma \beta(0)\sigma^{-1}=\beta(0).$ Thus $\sigma e^{\a}=de^{\a}$
for a nonzero constant $d\in \C.$ It is clear then that
$\sigma e^{n\alpha}=d^ne^{n\alpha}.$ Thus $\sigma=\sigma_d.$ This
completes the proof in the case $n=1.$ 

\begin{rem} {\rm In the case $n=2$ vertex operator algebra $(V_{L_4})^+$ is
isomorphic to $L(1/2,0)\otimes L(1/2,0)$ (see Lemma 3.1 of [DGH]).
One can easily see that the automorphism group of
$L(1/2,0)^{\otimes n}$ is the symmetric group $Sym_n$ for any $n>0.$ }
\end{rem} 

In the next section we will determine the automorphism group of $V_{L_8}^+.$

\section{$Aut(V_{L_8}^+)$} 

The proof in the case $n=4$ is more complicated and involves 
$(V_{L_2})_1$, which is a Lie algebra via the product $a_0b$.  
We denote it by $\frak g$; $\frak g \cong sl(2,C)$.  

We begin with a description of $Sym_4$ as a group of automorphisms of the
Lie algebra $\frak g$.  
We use the usual Chevalley basis $\a, x_\a, x_{-\a}$, which, as elements 
of $V_1 \le V$ are $\a(-1) \otimes e^0, 1\otimes e^\a, 1 \otimes e^{-\a}$.  

We define 
$\xap := x_\a \pm  x_{-\a}$ and $\xam := x_a - x_{\pm \a}$.  Then we
have 
$[\xap, \xam]=-2\a$ and for any scalar $c$, $[\a, c\xap]=2c\xam$ and 
$[\a, c\xam]=2c\xap$.  

Define $y_1 := \a, y_2 := \xap, y_3 := i\xam$.  Then for all $j$,
$[y_j,y_{j+1}]=y_{j+2}$ (indices read modulo 3).  

We have the double basis ${\cal D} :=   \{ \pm y_1, \pm y_2, \pm y_3 \}$.  
Its stabilizer in $Aut(\frak g)$ is a finite group $S \cong Sym_4$; 
the kernel in $S$ of its action on the three 1-spaces
spanned by elements of $\cal D$ is $E:=O_2(S)$, and $S$ induces the
full symmetric group on this set of three 1-spaces.  

The group $E$ consists of all evenly many sign changes at the three
1-spaces.  Define  $\t_j$ as the element of $E$ which takes $y_j$ to
itself and negates the other $y_k$.  Define, for each index $j$, the
two maps $\s_j$  and $\s'_j$ by 
$$\s_j : y_j \mt  -y_j; \ \ \ \ \ \   y_k \mt  y_\ell 
\mt y_k \hbox { if }  \{j,k,\ell\}=\{1,2,3\}.$$
and
$$\s_j' : y_j \mt  -y_j \hbox{ for } k=j; \ \ \ \ \ \   y_k \mt 
- -y_\ell \hbox{ and } y_\ell \mt  -y_k \hbox { if } 
\{j,k,\ell\}=\{1,2,3\}.$$
Then,  $\s_j'=\s_j \t_j = \t_j \s_j$, for all $j$.  Also, 
$\s_k \s_\ell$ maps  $y_k \mt y_\ell \mt -y_m$ and $\s_k \s_\ell
\s_k = \s_m'$, for $\{k,\ell, m\}=\{1,2,3\}$. 
It follows that $T:=\la \s_k, \s_\ell, \s_m' \ra \cong Sym_3$ and $T$
complements $E$ in $S$.  

\bigskip 

\def\am#1{\a (-#1)}
\def\ap#1{\a (#1)}
\def\etap{e^{2\a}+e^{-2\a}}
\def\etam{e^{2\a}-e^{-2\a}}
\def\eap{e^{\a}+e^{-\a}}
\def\eam{e^{\a}-e^{-\a}}
\def\peap{(e^{\a}+e^{-\a})}
\def\peam{(e^{\a}-e^{-\a})}

As pointed out in the introduction, $V_{L_8}^+ \cong V_{L_2}^{E}.$  So,
by Theorem 2.9,  we
need to determine the action of $S/E \cong Sym_4$ on $(V_{L_2}^E)_4.$ 
One can use Proposition 3.3 of [DM2] to see that different subgroups of $S/E$
give different subspaces of fixed points of $(V_{L_2}^E)_4.$ Thus
the action of $S/E $ on $(V_{L_2}^E)_4$ is faithful. Here 
we give a direct proof of this assertion.

\bigskip 
Let $V:=V_{L_8}^+$.  Then $dim(V_0), dim (V_1), \dots $ is $1,0,1,4,
\dots $.   The elements   $$\am1\am3, \ \am2^2,
\  \am1 ^4,   \ \etap$$    span $V_4$.   

\bigskip

Now to do some calculations.  We need some minimal information
about the action of certain
elements of $S$ on $V_4$.  
Take $\s := \s_2$, which negates $\eam$
and interchanges $\a$ and $\eap$.  So, $Y( \s ( \am1 ), z)\vv =Y(\eap,
z)\vv $.

\smallskip

We want to show that $\s_1$ and $\s$ act differently on $V_4$.
Clearly, $\s_1$ inverts $\a$ and so fixes  all the polynomial terms in the
above spanning set.  On the other hand, the 
coefficient of $z^4$ in $Y(\a, z)\vv$ (which is 
``all polynomial") goes under $\s$ to the 
coefficient of $z^4$ in $Y(\eap,z)\vv$.  
The latter coefficient has nontrivial contributions at the 
standard basis elements
$e^{2\a}$ 
and $e^{-2\a}$ of $V_{L_8}$, whence $\s$ does not leave the degree 4 polynomial
terms fixed, as $\s_1$ does.  This proves that these two 
involutions act differently on $V_4$.  We conclude that
$S/E \cong Sym_3$ acts faithfully on $V_4$.

\def\ea{e^\a}
\def\ema{e^{-\a}}  
\def\amf#1{  {{\am#1 \over #1}}  }
\def\amfz#1{  {{\am#1 \over #1}} z^#1  }
\def\amfp#1#2{  {{\am#1 ^#2   \over #1^#2 #2!}}
  }
\def\amfpz#1#2{  {{\am#1 ^#2  z^{#1 \cdot #2} \over #1^#2 #2!}}
  }

\bigskip

Next, we need a calculation to prove that
$(\am3\am1)_3(\am3\am1)=72\am3\am1$.

We need the coefficient at $z^{-4}$ of $Y(\am3\am1,z)(\am3\am1)$.  The
vertex operator is given by
$:[{1 \over 2 } {d^2 \over dz^2 } \a(z)] \a(z):$, where $\a(z) \colon
=  \sum_n \a(n) z^{-n-1}$ [FLM], (8.5.5).  We compute ${1\over 2}{d^2 \over dz^2
 }
\a(z) =\sum_n {n+2 \choose 2} \a(n) z^{-n-3}$, whence 
\begin{eqnarray*}
& &\ \ \ \ \ Y(\am3\am1,z)\\
& &=\sum_m \left[\sum_{i,j:i+j=m, i \le j} {i+2 \choose
2}\a(i)\a(j) +  \sum_{i,j:i+j=m, i > j} {i+2 \choose 2}\a(j)\a(i)\right]z^{-m-4}.
\end{eqnarray*}
When $m=i+j=0$ and the only summands which have nonzero value on 
$\am3\am1$ are those with $(i,j)=(-3,3),(-1,1),(1,-1),(3,-3)$.  Since
$(\a,\a)=2$ and ${-3 \choose 2}=1, {-1 \choose 2}=1, {1 \choose 2}=0$
and ${3 \choose 2}=2$, the standard  commutator relations 
$[h(m),h'(n)]=m(h,h')\delta_{m+n,0}$ imply that the coefficient of
$z^{-4}$ is $72 \am3\am1$.

\bigskip

The irreducible modules for $Sym_3$ are indicated by their dimensions
$1,1',2$, where $1$ is the trivial character and $1'$ is the sign
representation.  Tensor behavior includes these rules: $2 \otimes 2 =
2+1+1'$ and if the tensor of two irreducible contains $2$, then one
of the factors is $2$. 

\smallskip  

This information implies that the pairing $p:V_4 \times V_4 \rightarrow 
V_4$ given by the coefficient at $z^{-4}$ has the following property:
there is a decomposition $V_4=H \oplus J$, where $J$ is the
2-dimensional subspace spanned by $L(-2)^2, L(-4)$ and $H$ is the
unique nontrivial 
irreducible $S$-submodule in $V_4$; $J$ is a trivial module for
$Aut(V)$. 

\smallskip 
We argue that $H$ is invariant under $Aut(V)$.  First, $H$
contains a highest weight vector $\etap$ for the Virasoro element.   
{}From irreducibility of $S$ and its commuting with the Virasoro element,
$H$ consists of highest weight vectors. Note that  $J\subset L(1,0)$
does not contain any highest weight vector. 
Since $V_4=H \oplus J$, we immediately see that $H$ is the space of highest 
weight vectors, 
which clearly is invariant under $Aut(V).$ 

\smallskip

The pairing $p$ gives a nontrivial product on $H$ by 
$(x, y) \mapsto qp(x,y)$, where $q$ is the projection of $V_4=H \oplus
J$ to $H$.  Nontriviality of $qp$ follows from
$(\am3\am1)_3(\am3\am1)=72\am3\am1$ and the fact that 
$q(\am3\am1) \ne 0$ (because $\am3\am1$ is moved by $\s$)

\smallskip
 Now, we use the fact that the
(unique up to scalar) $Sym_3$-invariant algebra structure on $H$ has
automorphism group isomorphic to $Sym_3$ (see Appendix) to get that
$Aut(V) \cong Sym_3$, by restriction to $Aut(H,p)$.

\bigskip

\section{
Appendix. Noncommutative algebras of dimension $n-1$ with $Sym_n$ as
automorphism group. } 

\begin{thm}\label{cnaa}
If  $M$ is the $n-1$-dimensional 
irreducible part of the standard
degree $n$ permutation module for $Sym_n$, $M$ has a unique
(up to scalar multiplication) 
commutative algebra structure.  It is nonassociative and its
automorphism group is just $Sym_n$ by the given action.  
\end{thm} 

\smallskip 

The original proof of its automorphism group 
(RLG, 1977; unpublished preprint) amounted 
to showing that the set of idempotents $e \in M$ such
that $ad(e):x \mapsto xe$ has eigenvalues $1$ once and $-1 \over n-2$
with multiplicity $n-2$ is just an $n$-set which spans $M$. 

Here is a new proof.  
Let $(M,*)$ be the above algebra, invariant under $Sym_n$.  The
character theory of this irreducible and its low degree tensors says that 
an invariant nonzero commutative algebra structure is essentially unique.  
Here is one way to realize it.  Let $P$ be the 
permutation module on basis $e_i$ 
and make $P$ a direct sum of fields with indecomposable idempotents $e_i$.
This associative 
algebra has automorphism group isomorphic to 
$Sym_n$, seen by its action on the indecomposable idempotents.  
We get the submodule  $M$
as the annihilator of $1=\sum_i e_i$ with respect to the form with 
orthonormal basis $e_i$; this form is invariant under the automorphism
group since it may be realized as $(a,b) \mapsto tr(ad_P(a)ad_P(b))$.  
We get a commutative algebra structure on $M$ by using the given product on 
$P$, then orthogonally projecting the result to $M$.  This product in $M$ is 
nonzero for $n \ge 2$ 
since $e_1-e_2, e_1 - e_3 \in M$ has product 
in $P$ which is not a multiple of 
$1=\sum_i e_i$.  

The product on 
$P = M \oplus K.1$ is a map $P \otimes P \rightarrow P$ 
which may be expressed as a linear combination
of these maps:  
the natural identifications 
$K \otimes K = K$ and $K \otimes M = M$, the 
product $*$ on $M$, and the map  
$M \times M \rightarrow K$ 
which sends to $f(x,y)1$, where $f$ is the invariant inner
product.  
All such maps are invariant under $Aut(M,*)$, so $Aut(M,*)$
preserves the given algebra structure on 
$P$.  At once, we get $Aut(M,*) \cong Sym_n$.

\medskip
\noindent Department of Mathematics, University of California, 
Santa Cruz, CA 95064 USA. Email address: dong@cats.ucsc.edu (C.D.)

\noindent Department of Mathematics, University of Michigan, 
Ann Arbor, MI 48109-1109 USA. Email Address: rlg@math.lsa.umich.edu (R.L.G.)

\end{document}